\documentclass{amsart}
\usepackage{amsfonts}
\begin{document}

\newcommand{\C}{{\mathbb C}}
\newcommand{\G}{{\mathbb G}}
\newcommand{\GG}{{\check \mathbb G}}
\newcommand{\R}{{\mathbb R}}
\newcommand{\Q}{{\mathbb Q}}
\newcommand{\W}{{\mathbb W}}
\newcommand{\Z}{{\mathbb Z}}
\newcommand{\la}{{\langle}}
\newcommand{\ra}{{\rangle}}
\newcommand{\res}{{\rm res}}
\newcommand{\SP}{{\rm Sp}}
\newcommand{\half}{{\textstyle{\frac{1}{2}}}}
\newcommand{\tr}{{\rm Trace}}

\title {An algebraic analog of the Virasoro group} 

\author{Jack Morava}
\address{Department of Mathematics, Johns Hopkins University, Baltimore,
Maryland 21218}
\email{jack@math.jhu.edu}
\thanks{The author was supported in part by the NSF.}
\subjclass{81R10, 55S25}
\date {15 June 2001}

\begin{abstract} 
The group of diffeomorphisms of a circle is not
an infinite-dimensional algebraic group, though in many ways it 
behaves as if it were. Here we construct an algebraic model for this
object, and discuss some of its representations, which appear 
in the Kontsevich-Witten theory of two-dimensional topological 
gravity through the homotopy theory of moduli spaces. [This is 
a version of a talk on 23 June 2001 at the Prague Conference 
on Quantum Groups and Integrable Systems.] \end{abstract}

\maketitle

\section{Some functors from commutative rings to groups}

\noindent
{\bf 1.1} A formal diffeomorphism of the line, with coefficients in a 
commutative ring $A$, is an element $g$ of the ring $A[[x]]$ of formal 
power series with coefficients in $A$, such that $g(0) = 0$ and $g'(0)$ 
is a unit. More precisely, the {\bf group} of formal diffeomorphisms 
of the line, defined over $A$, is the set
\[
\G(A) = \{ g \in A[[x]] \;|\: g(x) = \sum_{k \geq 0} g_k x^{k+1} \; 
{\rm with} \; g_0 \in A^{\times} \; \} \;.
\]
Composition $g_0,g_1 \mapsto (g_0 \circ g_1)(x) = g_0(g_1(x))$ of formal 
power series makes this set into a monoid with $e(x) = x$ as identity 
element, and it is an exercise in induction to show that such an 
invertible power series [ie with leading coefficient a unit] possesses 
a composition inverse in $\G(A)$. Thus $\G$ defines a covariant functor from 
commutative rings to groups; in fact this functor is representable, in 
the sense that $\G(A)$ is naturally isomorphic to the set of ring 
homomorphisms from the polynomial algebra $\Z[g_k \;|\; k \geq 0][g_0^{-1}]$ 
to $A$. Composition endows this representing algebra with the Hopf diagonal 
\[
\Delta g(x) = (g \otimes 1)((1 \otimes g)(x)) \;,
\]
making $\G$ into a (pro-)algebraic group [9]. The kernel of the homomorphism 
\[
\epsilon \mapsto 0 : \G(A[\epsilon]/(\epsilon^2)) \to \G(A)
\]
can be given the structure of a Lie algebra, which is naturally isomorphic 
to the Lie algebra over $A$ spanned by the differentiation operators 
\[
v_k = x^{k+1}\frac{d}{dx}, \; k \geq 0 \;.
\]
satisfying $[v_k,v_l] = (l - k)v_{k+l}$. \bigskip

\noindent
{\bf 1.2} There is a closely related functor $\GG$ from commutative rings to 
groups, which in some ways resembles the group of diffeomorphisms of the 
circle. This functor is quite representable - it is an ind-proalgebraic 
group - but it is close enough to being so to have some useful properties. 
In particular its Lie algebra, in the sense above, is spanned by operators 
$v_k$ with $k \in \Z$: thus $k$ can be negative as well as positive. \bigskip

\noindent
In the terminology of [5 \S 2.3] an element $g$ of the Laurent series 
ring $A((x)) := A[[x]][x^{-1}]$ is a {\bf nil-Laurent} series if its 
coefficients $g_i$ for $i < -1$ are nilpotent. If $\surd A$ is the 
radical of $A$ (ie the ideal of nilpotent elements) and $A_{red} := 
A/\surd A$ then the set of such nil-Laurent series is the inverse image 
of $A_{red}[[x]]$ under the quotient homomorphism $\rho: A((x)) \to 
A_{red}((x))$, and 
\[
\GG(A) = \{ g \in A((x)) \:|\: \rho(g) \in \G(A_{red}) \}
\]
is the set of formal Laurent series $g(x) = \sum_{k \gg -\infty} g_k x^{k+1}
\in A((x))$ such that \medskip

i) $g_0$ is a unit, and 

ii) $g_{-k}$ is nilpotent, if $k \geq 1$. \medskip

\noindent
This set is closed under composition of power series, and is in fact
a group: We can write $g \in \GG$ as a sum 
\[
g(x) = g_+(x) + g_-(x^{-1})
\]
of an invertible formal power series $g_+$ and a polynomial $g_-$ in $x^{-1}$
with nilpotent coefficients; this implies that the sum
\[
g^{-1} = \sum_{k \geq 0} (-g_-)^k g_+^{-k-1} \in A((x))
\]
is finite. If $h = h_+ + h_-$ is another series of the same sort, we
can thus make sense of the composition $h_- \circ g$, so it suffices 
to show that $h_+ \circ (g_+ + g_-)$ is well-defined; but $g_-$, being 
a polynomial with nilpotent coefficients, is itself nilpotent, so this 
composition can be written as a finite Maclaurin expansion
\[
\sum_{k \geq 0} (D_k h_+)(g_+)g_-^k \in A((x)) \;,
\]
where
\[
D_k x^n = \frac{1}{k!}\frac{d^k x^n}{d x^k} = \binom{n}{k} x^{n-k} \; 
{\rm if} \; n \geq k \;, 
\]
and is otherwise zero. To show the existence of (composition) inverses, we 
use the fact that 
\[
h \circ g \equiv h \circ g_+ \; {\rm mod} \; I_-((x)) \;,
\]
where $I_-$ is the ideal generated by the coefficients of $g_-$. [I would
like to thank M. Kapranov for suggesting this line of argument, which has
substantially improved the result.] Because this ideal is generated by 
finitely many nilpotent elements, it is itself nilpotent, in the sense 
that $(I_-)^n = 0$ for $n \gg 0$. It suffices to construct an inverse for 
$g$ under the assumption that $g_0 = 1$, and that the rest of its 
coefficients lie in such a nil-ideal: for $g_+$ has a composition 
inverse $h^+$, such that 
\[
(g \circ h^+)(x) \equiv x \; {\rm mod} \; I_-((x)) \;,
\]
and if $u_0 \in A^{\times}$ is the coefficient of $x$ in $g \circ h^+$ 
then $h_{(0)}(x) = h^+(u_0^{-1}x)$ is a formal series such that $(g 
\circ h_{(0)})(x) \equiv x$ mod $I_-((x))$. Under that hypothesis, then,
let $h^+_{(1)}(x) = 2x - g(x) = x - \tilde g(x)$; then
\[
g(h^+_{(1)}(x)) = x - \tilde g(x) + \tilde g(x - \tilde g(x)) = x -
\tilde g'(x) \tilde g(x) + \dots \;,
\]
the dots representing further Taylor's series-style corrections,
so 
\[
g \circ h^+_{(1)} \equiv x \; {\rm mod} \; I_-^2((x)) \;.
\]
If $u_1$ is the coefficient of $x$ in $g \circ h^+_{(1)}$, and we
define $h_{(1)}(x) = h^+_{(1)}(u_1^{-1}x)$, then it follows that $g \circ 
h_{(1)}$ is a series with all coefficients in $I_-^2$ except that 
of $x$, which equals $1$. Now we can iterate this process: induction 
defines a sequence 
\[
h = h_{(0)} \circ h_{(1)} \circ \dots \;. 
\]
of compositions which will terminate in finitely many steps, defining 
the promised composition inverse for $g$. \bigskip

\noindent
{\bf 1.3} If $\phi$ is a diffeomorphism of the circle with the property
that 
\[
\phi(\zeta z) = \zeta \phi(z)
\]
(where $z \in \C$ with $|z| = 1$ and $\zeta$ is a primitive $n$th root
of unity), then $z \mapsto \phi(z)^n$ is an $n$-fold covering, which
factors through a diffeomorphism $\Phi$ of the circle satisfying
\[
\Phi(z^n) = \phi(z)^n \;.
\]
The group of such periodic diffeomorphisms thus defines an $n$-fold 
cover of Diff $S^1$. The group-valued functor $\GG$ has similar `covers':
for simplicity let $n=p$ be prime, e.g. two, and assume that $A$ contains 
a nontrivial $p$th root $\zeta$ of unity: then
\[
\GG_{1/p}(A) = \{ g \in \GG(A) \;|\; g(\zeta x) = \zeta g(x) \}
\]
is the subgroup of nil-Laurent series $g(x) = \sum_{k \gg -\infty} 
g_{kp} x^{kp+1}$ with $g_0$ a unit, and when $p$ is invertible in $A$
(e.g. if $A$ is a $\Q$-algebra) the homomorphism 
\[
g(x) \mapsto g(x^{1/p})^p := g^{(p)}(x) : \GG_{1/p} \to \GG
\]
induces an isomorphism of Lie algebras. This allows us to think of
the group of invertible nil-Laurent series in $A((x))$ as a subgroup of
the invertible nil-Laurent series in $A((x^{1/p}))$. \bigskip

\section{Some representations of these functors}

\noindent
Certain standard representations of Diff $S^1$ have analogs
for $\GG$; because these are representations over the complexes, I will
assume in this section that $A$ is an algebra over a field of characteristic 
zero. \bigskip

\noindent
{\bf 2.1} The $A$-bilinear form
\[
g,h \mapsto \la g,h \ra := - \sum_{k \in \Z} k \; g_{k-1} h_{-k-1} : A((x)) 
\times A((x)) \to A 
\]
is antisymmetric, and (aside from the subring of constants) is nondegenerate 
if $A$ is a $\Q$-algebra; it is an algebraic analog of the symplectic pairing 
\[
g,h \mapsto \res_{x=0} \; gdh
\]
of [10 \S 1]. The set $\SP_L(A)$ of `Laurent-symplectic' $A$-linear 
automorphisms of $A((x))$ which \medskip

i) preserve the bilinear form $\la \cdot,\cdot \ra$, and 

ii) are continuous in the pro-discrete topology of $A((x))$ \medskip

\noindent
defines a group-valued functor, analogous to the restricted symplectic group
[10 \S 5]. It is classical [12] that the residue of a differential over a 
local formal Laurent series field is independent of the choice of 
uniformizer. This remains true over general commutative rings $A$ [13], 
which implies that the composition 
\[
f \mapsto [h \mapsto h \circ f^{-1}] : \GG(A) \to \SP_L(A)
\]
is a natural homomorphism between group-valued functors; thus the 
$\GG$ has a natural linear representation as automorphisms of the functor 
which sends $A$ to the $A$-module $A((x))$. It is in any case elementary 
to see that the Lie algebra of $\GG$ preserves the symplectic form: if
\[
x \mapsto x + \epsilon x^{n+1}
\] 
then 
\[
x^k \mapsto x^k[1 + k \epsilon x^n] \;,\; dx^l \mapsto lx^{l-1}
[1 + (n+l)\epsilon x^n]dx
\]
so $\la x^k,x^l \ra$ changes under such a transformation by 
\[
l(n+k+l) \epsilon \; \res_{x=0} \; x^{n+k+l-1} dx \;.
\]
The residue in this expression is zero unless $n+k+l = 0$, 
but in that case the coefficient is zero; nothing in this 
argument requires that $n$ be positive. \bigskip

\noindent
{\bf 2.2} The residue pairing restricts to a bilinear form on the
ring $A[x,x^{-1}]$ of Laurent polynomials, which has a canonical
decomposition
\[
A[x,x^{-1}] = A[x] + A[x^{-1}]
\]
into Lagrangian subspaces. The symplectic form defines a Heisenberg algebra 
which is essentially (when $A$ is the field of real numbers) the identity
component of the loop group of the circle. The Fock representation [10 \S 3, 
11 \S 9.5] associated to this decomposition is an algebra of symmetric 
functions on the `positive-frequency' subspace $A[x]$. The restricted 
symplectic group acts as well on (a completion of) this representation, 
intertwining projectively with the action of the loops on the circle 
[10 \S 5, \S 7b; 11 \S 13.4]; this is simultaneously a (positive-energy) 
representation of the Heisenberg algebra of the bilinear form, and (an 
extension of) the Lie algebra of $\GG$. The Segal-Sugawara construction 
expresses the action of the Virasoro generators as quadratic expressions 
in the Heisenberg group elements [14 \S 1.7]. \bigskip

\noindent
{\bf 2.3} This Fock representation has an interpretation in terms of 
symmetric functions [4]; more generally, a certain class of twisted 
representations of the Virasoro algebra [2 \S 9.4], associated to 
Hall-Littlewood polynomials at roots of unity [7 III \S 8.12], fit 
naturally into this framework. For simplicity, let $p$ be a fixed prime 
[e.g. $p=2$] and let
\[
\C((x)) := V_0 \subset \C((x^{1/p})) := V
\]
be the extension of the field of formal complex Laurent series defined by 
adjoining a $p$th root of $x$. The Galois group of the field extension 
$V/V_0$ is cyclic of order $p$, generated by the automorphism
\[
x^{1/p} \mapsto \zeta x^{1/p} \;,
\]
and the bilinear form satisfies
\[
\la \zeta(g),\zeta(h) \ra = \la g,h \ra \;,
\]
so the invariant subspace $V_0$ is a bilinear submodule. More
generally, 
\[
V = \oplus_{a \in \Z/p\Z} V_a
\]
splits into orthogonal bilinear submodules $V_a$ spanned by
series of the form
\[
g = \sum_{s \gg -\infty} g_s x^s 
\]
in which $s = \pm(k + a/p)$ with $k$ and $a$ nonnegative
integers, $0 \leq a \leq p-1$. The restriction of the Fock representation 
of the group of nil-Laurent series over $\C((x^{1/p}))$ to $\GG_{1/p}$ 
can be interpreted (using the isomorphism of \S 1.3) as a representation 
of $\GG$ on the (completed) tensor product of the rings $S(V_a)$ of 
symmetric functions. \bigskip

\noindent
These rings have Hopf algebra structures, which are usually described 
in terms of exponentials: the Witt functor assigns to a commutative ring 
$A$ the multiplicative group 
\[
\W(A) = (1 + xA[[x]])^{\times}
\]
of formal series $w(x) = \sum w_k x^k$ with constant coefficient
$w_0 = 1$. This is naturally isomorphic to the set of ring
homomorphisms from a polynomial algebra on generators $\{w_k,\; k
> 0\}$ to $A$; the group structure endows this representing ring
with the structure of a (commutative and cocommutative) Hopf
algebra. The involution
\[
w(T) \mapsto w^*(x) = w(-x)
\]
respects the product, and so defines a $\Z/2\Z$-action on $\W$.
The Hopf algebra of Schur $Q$-functions represents the kernel of the
norm homomorphism
\[
w \mapsto  w \cdot w^* : \W \to \W \;;
\]
in other words it represents the functor which sends a ring to
the group of power series $q(x)$ with $q(0) = 1$ over that ring,
which satisfy the relation $q(x)q(-x) = 1$. \bigskip

\noindent
This ring is torsion-free, and we can reformulate the
relation above in the universal example as the assertion that
the formal logarithm $\log q(x)$ is an odd power series in $T$.
More generally, the group of $p$th roots of unity acts on $\W$ by
$w(x) \mapsto w(\zeta x)$, and 
\[
w(x) \mapsto \prod_{a \in \Z/pZ} w(\zeta_p^a x) = N(w)(x) :\W \to
\W 
\]
is the Frobenius homomorphism of Witt theory. The Hopf algebra
representing its kernel can be described as an algebra of Hall-Littlewood
symmetric functions evaluated at a $p$th root of unity. For our purposes 
it can most conveniently be understood in terms of power series $w(x)$ with 
$w(0)=1$ such that the projection of $\log w(x^{1/p})$ to $V_0$ is zero. 
The polynomial algebra underlying the Fock representation thus splits 
(over $\Q$) as a product of Hopf algebras, indexed by $a \in 
\Z/p\Z$; its $a$th component is the Fock representation of 
the Heisenberg group defined by $V_a$. \bigskip

\noindent
{\bf 2.4} The primitives in these Hopf algebras acquire natural normalizations
from the Heisenberg algebra: the fractional divided powers
\[
\gamma_s = \frac{x^s}{\Gamma(s+1)}
\]
satisfy
\[
\la \gamma_{n/p},\gamma_{m/p} \ra = \frac{\res \; x^{(n+m)/p-1} \;
dx}{\Gamma(1+n/p)\Gamma(m/p)} = - \pi^{-1} \delta_{n+m,0} \; \sin
(n\pi/p) \;,   2 
\] 
so the elements
\[
\gamma(a)_{\pm k} := |\sin(a\pi/p)|^{-\half}
\gamma_{\pm(k+a/p)} \;, k \in \Z_+
\]
define a normalized symplectic basis for $V_a$ when $a$ is {\bf not}
congruent to zero mod $p$. \bigskip

\noindent
When $p=2$, this defines the Virasoro representation with $c=1$ and 
$h=1/16$ studied in the Kontsevich-Witten theory of two-dimensional 
topological gravity [1]; that theory has a conjectural generalization 
[3,6] in which the more general representations defined above (with 
$c=1$ and $h = (p^2-1)/48$) play a similar role. From a geometric point 
of view, these representations are mysterious: they are somehow 
homotopy-theoretic, and do not arise in any natural way from the Lie 
algebra of vector fields on the circle; instead, they seem to be 
related to automorphisms of the cohomology of infinite-dimensional 
complex projective space, along the lines laid out in this paper, t
hrough work of Madsen and Tillmann [8].

\bibliographystyle{amsplain}

\end{document}